\documentclass{article}

\newtheorem{theorem}{Theorem}[section]
\newtheorem{corollary}[theorem]{Corollary}

\newtheorem{lemma}[theorem]{Lemma}
\newtheorem{proposition}[theorem]{Proposition}

\evensidemargin\oddsidemargin
\begin{document}

\author{Vadim E. Levit and Eugen Mandrescu \\
%EndAName
Department of Computer Science\\
Holon Academic Institute of Technology\\
52 Golomb Str., P.O. Box 305\\
Holon 58102, ISRAEL\\
\{levitv, eugen\_m\}@barley.cteh.ac.il}
\title{On $\alpha ^{++}$-Stable Graphs}
\maketitle

\begin{abstract}
The \textit{stability number }of the graph $G$, denoted by $\alpha (G)$, is
the cardinality of a maximum stable set of $G$. A graph is well-covered if
every maximal stable set has the same size. $G$ is a \textit{%
K\"{o}nig-Egerv\'{a}ry graph} if its order equals $\alpha (G)+\mu (G)$,
where $\mu (G)$ is the cardinality of a maximum matching in $G$. In this
paper we characterize $\alpha ^{++}$-stable graphs, namely, the graphs whose
stability numbers are invariant to adding any two edges from their
complements. We show that a K\"{o}nig-Egerv\'{a}ry graph is $\alpha ^{++}$%
-stable if and only if it has a perfect matching consisting of pendant edges
and no four vertices of the graph span a cycle. As a corollary it gives
necessary and sufficient conditions for $\alpha ^{++}$-stability of
bipartite graphs and trees. For instance, we prove that a bipartite graph is 
$\alpha ^{++}$-stable if and only if it is well-covered and $C_{4}$-free.
\end{abstract}

\section{Introduction}

All the graphs considered in this paper have at least two vertices. For such
a graph $G=(V,E)$ we denote its vertex set by $V=V(G)$ and its edge set by $%
E=E(G).$ If $X\subset V$, then $G[X]$ is the subgraph of $G$ spanned by $X$.
By $G-W$ we mean the subgraph $G[V-W]$ , if $W\subset V(G)$. By $G-F$ we
denote the partial subgraph of $G$ obtained by deleting the edges of $F$,
for $F\subset E(G)$, and we use $G-e$, if $W$ $=\{e\}$. If $A,B$ $\subset V$
and $A\cap B=\emptyset $, then $(A,B)$ stands for the set $\{e=ab:a\in
A,b\in B,e\in E\}$. The neighborhood of a vertex $v\in V$ is the set $%
N(v)=\{w:w\in V$ \ \textit{and} $vw\in E\}$, and $N(A)=\cup \{N(v):v\in A\}$%
, for $A\subset V$. If $\left| N(v)\right| =\left| \{w\}\right| =1$, then $v$
is a \textit{pendant vertex} and $vw$ is a \textit{pendant edge} of $G$. By $%
C_{n}$, $K_{n}$, $P_{n}$ we denote the chordless cycle on $n\geq $ $4$
vertices, the complete graph on $n\geq 2$ vertices, and respectively the
chordless path on $n\geq 3$ vertices.

A stable set of maximum size will be referred as to a \textit{maximum stable
set} of $G$. The \textit{stability number }of $G$, denoted by $\alpha (G)$,
is the cardinality of a maximum stable set in $G$. Let $\Omega (G)$ denotes $%
\{S:S$ \textit{is a maximum stable set of} $G\}$ and $\xi (G)=\left| \cap
\{S:S\in \Omega (G)\}\right| $. We call $\{\Omega _{1},\Omega _{2}\}$ a
cover of $\Omega (G)$ if $\Omega _{1},\Omega _{2}\subset \Omega (G)$ and $%
\Omega _{1}\cup \Omega _{2}=\Omega (G)$; by $\xi (\Omega _{i}),i=1,2$, we
mean the number $\left| \cap \{S:S\in \Omega (G_{i})\}\right| $.

A \textit{matching} is a set of non-incident edges of $G$; a matching of
maximum cardinality $\mu (G)$ is a \textit{maximum matching}, and a \textit{%
perfect matching} is a matching covering all the vertices of $G$. $G$ is a 
\textit{K\"{o}nig-Egerv\'{a}ry graph }provided $\alpha (G)+\mu (G)=\left|
V(G)\right| $, \cite{dem}, \cite{ster}.

A graph $G$ is $\alpha ^{+}$-\textit{stable} if $\alpha (G+e)=\alpha (G)$,
for any edge $e\in E(\overline{G})$, where $\overline{G}$ is the complement
of $G$, \cite{gun}. Haynes et al. have characterized the $\alpha ^{+}$%
-stable as follows:

\begin{theorem}
\label{th2}\cite{hayn} A graph $G$ is $\alpha ^{+}$-stable if and only if $%
\xi (G)\leq 1$.
\end{theorem}

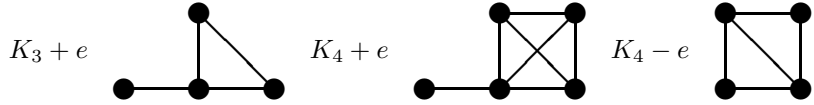
\begin{figure}[h]
\setlength{\unitlength}{1cm}%
\begin{picture}(5,1.2)\thicklines
  \multiput(2,0)(1,0){3}{\circle*{0.29}}
  \put(3,1){\circle*{0.29}}
  \put(2,0){\line(1,0){2}}
  \put(3,1){\line(1,-1){1}}
  \put(3,0){\line(0,1){1}} 
 \put(1,0.5){\makebox(0,0){$K_{3}+e$}} 

  \multiput(6,0)(1,0){3}{\circle*{0.29}}
  \multiput(7,1)(1,0){2}{\circle*{0.29}}
  \put(6,0){\line(1,0){2}}
  \put(7,1){\line(1,0){1}}
  \put(7,1){\line(1,-1){1}} 
  \put(7,0){\line(1,1){1}}
  \put(7,0){\line(0,1){1}}
  \put(8,0){\line(0,1){1}}

\put(5,0.5){\makebox(0,0){$K_{4}+e$}}

  \multiput(10,0)(1,0){2}{\circle*{0.29}}
  \multiput(10,1)(1,0){2}{\circle*{0.29}}
  \put(10,0){\line(1,0){1}}
  \put(10,1){\line(1,0){1}}
  \put(10,1){\line(1,-1){1}} 
  \put(10,0){\line(0,1){1}}
  \put(11,0){\line(0,1){1}}

\put(9,0.5){\makebox(0,0){$K_{4}-e$}}

 \end{picture}
\caption{$K_{3}+e$ and $K_{4}+e$ are $\alpha _{1}^{+}$-stable graphs; $%
K_{4}-e$ is not $\alpha ^{+}$-stable.}
\label{fig2}
\end{figure}

Theorem \ref{th2} implies that for an $\alpha ^{+}$-stable graph either $\xi
(G)=0$ or $\xi (G)=1$. This motivates the following definition. A graph $G$
is called ($\mathit{i}$) $\alpha _{0}^{+}$-\textit{stable} whenever $\xi
(G)=0$, and ($\mathit{ii}$) $\alpha _{1}^{+}$-\textit{stable} provided $\xi
(G)=1$, \cite{levm1}. For instance, $C_{4}$ is $\alpha _{0}^{+}$-stable,
while the graphs $K_{3}+e,K_{4}+e$ in Figure \ref{fig2} are $\alpha _{1}^{+}$%
-stable.

In \cite{hayn} it was shown that an $\alpha ^{+}$-stable tree can be only $%
\alpha _{0}^{+}$-stable, and this is exactly the case of trees possessing a
perfect matching. This result was generalized to bipartite graphs in \cite
{levm1}. Nevertheless, there exist both $\alpha _{1}^{+}$-stable
K\"{o}nig-Egerv\'{a}ry graphs (e.g., the graph $K_{3}+e$ in Figure \ref{fig2}%
), and $\alpha _{0}^{+}$-stable K\"{o}nig-Egerv\'{a}ry graphs (e.g., all $%
\alpha ^{+}$-stable bipartite graphs). A necessary (but not sufficient,
e.g., $K_{4}-e$) condition for $\alpha ^{+}$-stability is:

\begin{proposition}
\label{prop1}\cite{levm2} Any $\alpha ^{+}$-stable K\"{o}nig-Egerv\'{a}ry
graph has a perfect matching.
\end{proposition}

Let define a graph $G$ as $\alpha ^{++}$\textit{-stable} whenever $\alpha
(G+e_{1}+e_{2})=\alpha (G)$, for any $e_{1},e_{2}\in E(\overline{G})$, and $%
\alpha _{P_{3}}^{+}$-\textit{stable} provided $\alpha (G+e_{1}+e_{2})=\alpha
(G)$, for any $e_{1},e_{2}\in E(\overline{G})$ that have a common endpoint.
Gunther et al., \cite{gun}, studied the structure of $\alpha ^{+}$-stable
trees, and in \cite{levm1}, \cite{levm2}, some of their results were
extended to bipartite graphs and K\"{o}nig-Egerv\'{a}ry graphs.

In this paper we characterize $\alpha ^{++}$-stable graphs. We settle a
number of connections between $\alpha ^{+}$-stable graphs, $\alpha
_{P_{3}}^{+}$-stable and $\alpha ^{++}$-stable graphs. In particular, we
show that any $\alpha _{0}^{+}$-stable graph is $\alpha _{P_{3}}^{+}$%
-stable. We also give a necessary and sufficient condition for a graph to be 
$\alpha _{1}^{+}$-stable and $\alpha _{P_{3}}^{+}$-stable at the same time.

We prove that a K\"{o}nig-Egerv\'{a}ry graph is $\alpha ^{++}$-stable if and
only if it has a perfect matching consisting of only pendant edges and
contains no cycle on $4$ vertices. Using this result we describe $\alpha
^{++}$-stable bipartite graphs and $\alpha ^{++}$-stable trees. For
instance, it is shown that a bipartite graph is $\alpha ^{++}$-stable if and
only if it is well-covered and $C_{4}$-free.

\section{$\alpha ^{++}$-stable graphs}

Notice that any $\alpha ^{++}$-stable graph is $\alpha _{P_{3}}^{+}$-stable,
but the converse is not generally true. For instance, $C_{4}$ is $\alpha
_{P_{3}}^{+}$-stable but not $\alpha ^{++}$-stable. Let us also observe that
for $n\geq 2,$ the graph $K_{n}-e$ is $\alpha _{P_{3}}^{+}$-stable, but it
is not $\alpha ^{+}$-stable.

\begin{proposition}
\label{prop10}If $G\neq K_{n}-e,n\geq 2$ is $\alpha _{P_{3}}^{+}$-stable,
then it is also $\alpha ^{+}$-stable.
\end{proposition}

\setlength {\parindent}{0.0cm}\textbf{Proof.} Assume that $\alpha (G)=2$.
Since $G\neq K_{n}-e$, it follows that $\left| \Omega (G)\right| \geq 2$,
and consequently $G$ is $\alpha ^{+}$-stable, as well. For $\alpha (G)\geq 3$%
, suppose, on the contrary, that $G$ is not $\alpha ^{+}$-stable, and let $%
x,y\in \cap \{S:S\in \Omega (G)\}$. Hence, for $S\in \Omega (G)$ and $z\in
S-\{x,y\}$, we obtain that $\alpha (G+xy+xz)<\alpha (G)$, in contradiction
with the fact that $G$ is $\alpha _{P_{3}}^{+}$-stable. Therefore, $G$ must
be $\alpha ^{+}$-stable. \rule{2mm}{2mm}\setlength
{\parindent}{3.45ex}\newline

It is worth observing that an $\alpha ^{+}$-stable graph is not necessarily $%
\alpha _{P_{3}}^{+}$-stable. For instance, $K_{3}+e$ is $\alpha ^{+}$%
-stable, in fact it is $\alpha _{1}^{+}$-stable, but it is not $\alpha
_{P_{3}}^{+}$-stable. However, there exist graphs that are both $\alpha
_{1}^{+}$-stable and $\alpha _{P_{3}}^{+}$-stable; e.g., the graph $K_{4}+e$
in Figure \ref{fig2}.

\begin{proposition}
\label{prop8}Any $\alpha _{0}^{+}$-stable graph is also $\alpha _{P_{3}}^{+}$%
-stable.
\end{proposition}

\setlength {\parindent}{0.0cm}\textbf{Proof.} Suppose, on the contrary, that
some $\alpha _{0}^{+}$-stable graph $G$ is not $\alpha _{P_{3}}^{+}$-stable,
i.e., there are $x,y,z\in V$ such that $\alpha (G+xy+xz)<\alpha (G)$. Hence,
it follows that $x\in \cap \{S:S\in \Omega (G)\}$, in contradiction with the
fact that $G$ is $\alpha _{0}^{+}$-stable. \rule{2mm}{2mm}%
\setlength
{\parindent}{3.45ex}\newline

Combining Theorem \ref{th2} and Propositions \ref{prop10}, \ref{prop8}, we
get:

\begin{corollary}
If $G\neq K_{n}-e$ has $\alpha (G)=2$, then $G$ is:

($\mathit{i}$) $\alpha ^{+}$-stable if and only if $\left| \Omega (G)\right|
\geq 2$;

($\mathit{ii}$) $\alpha _{P_{3}}^{+}$-stable if and only if either it is $%
\alpha _{0}^{+}$-stable or it is $\alpha _{1}^{+}$-stable and $\left| \Omega
(G)\right| \geq 3$;

($\mathit{iii}$) $\alpha ^{++}$-stable if and only if $\left| \Omega
(G)\right| \geq 3$;

($\mathit{iv}$) $\alpha _{P_{3}}^{+}$-stable if and only if it is $\alpha
^{+}$-stable and $G\neq P_{3}(K_{1},K_{m},K_{2})$, where $%
P_{3}(K_{1},K_{m},K_{2})$ is the graph obtained by substituting the vertices
of $P_{3}$ respectively, by $K_{1},K_{m},K_{2}$, and joining all the
vertices of $K_{m}$ with the two vertices of $K_{2}$ and the single vertex
of $K_{1}$.
\end{corollary}

\begin{proposition}
\label{prop2}If $G=(V,E)$ has $\alpha (G)\geq 3$, then the following
assertions are equivalent:

\textit{(}$\mathit{i}$\textit{)} $G$ is $\alpha ^{+}$-stable;

($\mathit{ii}$) either $G$ is $\alpha _{P_{3}}^{+}$-stable, or there exist
three vertices $x,y,z\in V$ such that $\left| \{x,y\}\cap S\right| \cdot
\left| \{x,z\}\cap S\right| \geq 2$ holds for any $S\in \Omega (G)$, and $x$
is the unique vertex of $G$ with this property.
\end{proposition}

\setlength {\parindent}{0.0cm}\textbf{Proof.} ($\mathit{i}$) $\Rightarrow $ (%
$\mathit{ii}$) Let $G$ be $\alpha ^{+}$-stable, but not $\alpha _{P_{3}}^{+}$%
-stable. Hence, there are $x,y,z\in V$ such that $\alpha (G+xy+xz)<\alpha
(G) $. Therefore, we get that $x\in \cap \{S:S\in \Omega (G)\}$, because
otherwise, any $S\in \Omega (G)$ not containing $x$ is still stable in $%
G+xy+xz$, and consequently, we obtain $\alpha (G+xy+xz)=\alpha (G)$, in
contradiction with the assumption on $G$. In addition, each $S\in \Omega (G)$
satisfies $\left| S\cap \{y,z\}\right| \geq 1$, since otherwise, if some $%
S_{0}\in \Omega (G)$ has $S_{0}\cap \{y,z\}=\emptyset $, then $S_{0}$ is
stable in $G+xy+xz$ and this yields $\alpha (G+xy+xz)=\alpha (G)$, again in
contradiction with the assumption on $G$. Finally, $x$ is unique, because
otherwise $\xi (G)\geq 2$, which contradicts the $\alpha ^{+}$-stability of $%
G$.\setlength
{\parindent}{3.45ex}

($\mathit{ii}$) $\Rightarrow $ ($\mathit{i}$) If $G$ is $\alpha _{P_{3}}^{+}$%
-stable and $\alpha (G)\geq 3$, then by Proposition \ref{prop10}, $G$ is $%
\alpha ^{+}$-stable. Further, if there are $x,y,z\in V$ such that $\left|
\{x,y\}\cap S\right| \cdot \left| \{x,z\}\cap S\right| \geq 2$ holds for any 
$S\in \Omega (G)$, and $x$ is unique with this property, it follows that $G$
is not $\alpha _{P_{3}}^{+}$-stable, because $\alpha (G+xy+xz)<\alpha (G)$,
but it is $\alpha _{1}^{+}$-stable, since $\{x\}=\cap \{S:S\in \Omega (G)\}$%
. \rule{2mm}{2mm}

\begin{lemma}
\label{lem3}If for any $x,y\in V(G)$ there exists $S\in \Omega (G)$ such
that $x,y\in V(G)-S$, then $G$ is both $\alpha _{0}^{+}$-stable and $\alpha
^{++}$-stable.
\end{lemma}

\setlength {\parindent}{0.0cm}\textbf{Proof.} Suppose, on the contrary, that
there exists $x\in \cap \{S:S\in \Omega (G)\}$. Then, for any $y\in
V(G)-\{x\}$ and $S\in \Omega (G)$, we get $\{x,y\}\cap S\neq \emptyset $, in
contradiction with the premises on $G$. Therefore, $G$ is $\alpha _{0}^{+}$%
-stable.\setlength
{\parindent}{3.45ex} According to Proposition \ref{prop8}, it follows that $%
G $ is $\alpha _{P_{3}}^{+}$-stable, too. Assume that $G$ is not $\alpha
^{++}$-stable. Hence, since $G$ is $\alpha _{P_{3}}^{+}$-stable, we infer
that there are $x,y,u,v\in V(G)$, pairwise distinct, such that $\alpha
(G+xy+uv)<\alpha (G)$. Let $S\in \Omega (G)$ be such that $x,u\in V(G)-S$.
Then $S$ is stable in $G+xy+uv$, in contradiction with the assumption on $G$%
. Consequently, $G$ is $\alpha ^{++}$-stable. \rule{2mm}{2mm}\newline

As an example, $C_{2k+1},k\geq 2$ and $K_{n},n\geq 3$ are both $\alpha
_{0}^{+}$-stable and $\alpha ^{++}$-stable, according to Lemma \ref{lem3}.
Notice that the converse of Lemma \ref{lem3} is not generally true; see, for
instance, the graphs $C_{2k},k\geq 3$. There exist $\alpha _{0}^{+}$-stable
graphs that are not $\alpha ^{++}$-stable (e.g., $C_{4}$), and vice-versa,
there are $\alpha ^{++}$-stable that are not $\alpha _{0}^{+}$-stable (e.g., 
$K_{4}+e$).

\begin{proposition}
Let $G$ be $\alpha _{1}^{+}$-stable, $\{v\}=\cap \{S:S\in \Omega (G)\}$ and $%
G_{0}=G-N[v]$. If $G$ is not $\alpha _{P_{3}}^{+}$-stable, then there are $x$
and $y$ belonging to the same connected component of $G_{0}$, such that $%
\alpha (G+xv+yv)<\alpha (G)$. In other words, there exists a path connecting 
$x$ and $y$, which avoid the neighborhood of $v$.
\end{proposition}

\setlength {\parindent}{0.0cm}\textbf{Proof.} Let $\{H_{k}:1\leq k\leq
s\},s\geq 2$, be the connected components of $G_{0}$, and suppose that there
are $x$ and $y$ belonging to different connected components of $G_{0}$, (say
respectively $H_{i},H_{j}$), such that $\alpha (G+xv+yv)<\alpha (G)$. Since $%
G_{0}$ is $\alpha _{0}^{+}$-stable, it follows that all its connected
components are $\alpha _{0}^{+}$-stable, as well. Let $S_{k}$ $\in \Omega
(H_{k}),1\leq k\leq s$, and $S_{i}$ $\in \Omega (H_{i}),S_{j}$ $\in \Omega
(H_{j})$ be such that $x\notin S_{i},y\notin S_{j}$, which exist, because
all $H_{k}$ are $\alpha _{0}^{+}$-stable. Hence, we get that $\{v\}\cup
(\cup \{S_{k}:1\leq k\leq s\})$ is stable in $G+xv+yv$, in contradiction
with $\alpha (G+xv+yv)<\alpha (G)$. \rule{2mm}{2mm}%
\setlength
{\parindent}{3.45ex}

\begin{theorem}
Let $G$ be $\alpha _{1}^{+}$-stable, $\{v\}=\cap \{S:S\in \Omega (G)\}$ and $%
G_{0}=G-N[v]$. Then $G$ is $\alpha _{P_{3}}^{+}$-stable if and only if for
every pair $x,y\in V(G_{0})$ there exists $S_{0}\in \Omega (G_{0})$ such
that $x,y\in V(G_{0})-S_{0}$.
\end{theorem}

\setlength {\parindent}{0.0cm}\textbf{Proof.} Let $x,y\in V(G_{0})$. By
definition of $G_{0}$, it follows that $x,y\notin N[v]$, and since $G$ is $%
\alpha _{P_{3}}^{+}$-stable, we get that $\alpha (G+xv+yv)=\alpha (G)$.
Therefore, there is $S\in \Omega (G)$ such that $x,y\in V(G)-S$. Hence, $%
x,y\in V(G_{0})-S_{0}$, where $S_{0}=S-\{v\}\in \Omega (G_{0})$.%
\setlength
{\parindent}{3.45ex}

Conversely, $G$ is $\alpha _{1}^{+}$-stable, and for every pair $x,y\in
V(G_{0})$ there exists $S\in \Omega (G_{0})$ with $x,y\in V(G_{0})-S_{0}$.
Assume that $G$ is not $\alpha _{P_{3}}^{+}$-stable. Hence, there are $%
x,y\in V(G)$ such that $\alpha (G+xv+yv)<\alpha (G)$, since $G$ is $\alpha
_{1}^{+}$-stable. Let $S_{0}\in \Omega (G_{0})$ be such that $x,y\in
V(G_{0})-S_{0}$. Then, it follows that $S_{0}\cup \{v\}\in \Omega (G)\}$, in
contradiction with the assumption on $G$. Therefore, $G$ is also $\alpha
_{P_{3}}^{+}$-stable. \rule{2mm}{2mm}

\begin{proposition}
\label{prop3}A graph $G$ is not $\alpha _{P_{3}}^{+}$-stable if and only if $%
\xi (G)\geq 1$ and there exists a cover $\{\Omega _{1},\Omega _{2}\}$ of $%
\Omega (G)$, such that $\xi (\Omega _{i})\geq 2,i=1,2$.
\end{proposition}

\setlength {\parindent}{0.0cm}\textbf{Proof.} If $G$ is not $\alpha
_{P_{3}}^{+}$-stable, then $\alpha (G+e_{1}+e_{2})<\alpha (G)$ holds for
some $e_{1},e_{2}\in E(\overline{G})$ that have a common endpoint. Suppose $%
e_{1}=xy,e_{2}=yz$. Let us define 
\[
\Omega _{1}=\{S:x,y\in S\in \Omega (G)\}\ and\ \Omega _{2}=\{S:y,z\in S\in
\Omega (G)\}. 
\]
Hence, it follows that $\xi (G)\geq 1$ and $\xi (\Omega _{i})\geq 2,i=1,2$.%
\setlength
{\parindent}{3.45ex}

Conversely, assume that $\xi (G)\geq 1$, i.e., there exists at least one
vertex belonging to $\cap \{S:S\in \Omega (G)\}$, say $y$, and that there is
some cover $\{\Omega _{1},\Omega _{2}\}$ of $\Omega (G)$, such that $\xi
(\Omega _{i})\geq 2,i=1,2$. If $x\in \cap \{S:S\in \Omega _{1}\}-\{y\}$ and $%
v\in \cap \{S:S\in \Omega _{2}\}-\{y\}$, then $\alpha (G+xy+uv)<\alpha (G)$,
because any $S\in \Omega (G)$ contains at least one of the pairs $\{x,y\}$
or $\{y,v\}$. Therefore, $G$ can not be $\alpha _{P_{3}}^{+}$-stable. \rule%
{2mm}{2mm}

\begin{proposition}
\label{prop4}A graph $G$ is not $\alpha ^{++}$-stable if and only if there
exists a cover $\{\Omega _{1},\Omega _{2}\}$ of $\Omega (G)$, such that $\xi
(\Omega _{i})\geq 2,i=1,2$.
\end{proposition}

\setlength {\parindent}{0.0cm}\textbf{Proof.} If $G$ is not an $\alpha ^{++}$%
-stable graph, then $\alpha (G+e_{1}+e_{2})<\alpha (G)$ holds for some $%
e_{1},e_{2}\in E(\overline{G})$. Suppose $e_{1}=xy,e_{2}=uv$. Let us define $%
\Omega _{1}=\{S:x,y\in S\in \Omega (G)\}$ and $\Omega _{2}=\{S:u,v\in S\in
\Omega (G)\}$. Hence, it follows that $\xi (\Omega _{i})\geq 2,i=1,2$.
Suppose that there exists $S\in \Omega -\left( \Omega _{1}\cup \Omega
_{2}\right) $. Then $S\in \Omega \left( G+e_{1}+e_{2}\right) $, that
contradicts the inequality $\alpha (G+e_{1}+e_{2})<\alpha (G)$. Hence, $%
\Omega _{1}\cup \Omega _{2}=\Omega $, which means that $\{\Omega _{1},\Omega
_{2}\}$ is a cover we were supposed to find.\setlength
{\parindent}{3.45ex}

Conversely, assume that there is a cover $\{\Omega _{1},\Omega _{2}\}$ of $%
\Omega (G)$ with $\xi (\Omega _{i})\geq 2,i=1,2$. If $x,y\in \cap \{S:S\in
\Omega _{1}\}$ and $u,v\in \cap \{S:S\in \Omega _{2}\}$, then $\alpha
(G+xy+uv)<\alpha (G)$, because any $S\in \Omega (G)$ contains at least one
of the pairs $\{x,y\}$ or $\{u,v\}$. Therefore, $G$ is not $\alpha ^{++}$%
-stable. \rule{2mm}{2mm}\newline

Combining Propositions \ref{prop2} and \ref{prop4}, we deduce the following:

\begin{theorem}
For a graph $G$ the following assertions are equivalent:

($\mathit{i}$) $G$ is $\alpha ^{++}$-stable;

($\mathit{ii}$) $G$ is $\alpha ^{+}$-stable and $\Omega (G+e_{1})\cap \Omega
(G+e_{1})\neq \emptyset $ for any $e_{1},e_{2}\in E(\overline{G})$;

($\mathit{iii}$) $\Omega (G)\cap \Omega (G+e_{1})\cap \Omega (G+e_{1})\neq
\emptyset $ for any $e_{1},e_{2}\in E(\overline{G})$;

($\mathit{iv}$) $G$ is $\alpha ^{+}$-stable and $\left| \cap \{S:S\in \Omega
(G+e)\right| \leq 1$ for any $e\in E(\overline{G})$;

($\mathit{v}$) $G$ is $\alpha _{P_{3}}^{+}$-stable and there are no $%
e_{1},e_{2}\in E(\overline{G}),e_{1}=xy,e_{2}=uv$, such that: 
\[
\{x,y\}\cap \{u,v\}=\emptyset ,\ and\ for\ any\ S\in \Omega (G),\max
\{\left| S\cap \{x,y\}\right| ,\left| S\cap \{u,v\}\right| \}=2; 
\]

($\mathit{vi}$) for any cover $\{\Omega _{1},\Omega _{2}\}$ of $\Omega (G)$
either $\xi (\Omega _{1})\leq 1$ or $\xi (\Omega _{2})\leq 1$ holds.
\end{theorem}

\section{$\alpha ^{++}$-stable K\"{o}nig-Egerv\'{a}ry graphs}

According to a well-known result of K\"{o}nig, \cite{koen}, and
Egerv\'{a}ry, \cite{eger}, any bipartite graph is a K\"{o}nig-Egerv\'{a}ry
graph. This class includes also non-bipartite graphs (see, for instance, the
graph $K_{3}+e$ in Figure \ref{fig2}). If $G_{i}=(V_{i},E_{i}),i=1,2$, are
two disjoint graphs, then $G=G_{1}*G_{2}$ is defined as the graph with $%
V(G)=V(G_{1})\cup V(G_{2})$, and 
\[
E(G)=E(G_{1})\cup E(G_{2})\cup \{xy:for\ some\ x\in V(G_{1})\ and\ y\in
V(G_{2})\}. 
\]

\begin{proposition}
\cite{levm2}\label{prop11} The following assertions are equivalent:

($\mathit{i}$) $G$ is a K\"{o}nig-Egerv\'{a}ry graph;

($\mathit{ii}$) $G=H_{1}*H_{2}$, where $V(H_{1})=S\in \Omega (G)$ and $%
\left| V(H_{1})\right| \geq \mu (G)=\left| V(H_{2})\right| $;

($\mathit{iii}$) $G=H_{1}*H_{2}$, where $V(H_{1})=S$ is a stable set in $G,$ 
$\left| S\right| \geq \left| V(H_{2})\right| $ and $(S,V(H_{2}))$ contains a
matching $M$ with $\left| M\right| =\left| V(H_{2})\right| $.
\end{proposition}

It it easy to see that a K\"{o}nig-Egerv\'{a}ry graph $G$ has a perfect
matching if and only if $\alpha (G)=\mu (G)$. The edges of any maximum
matching of a K\"{o}nig-Egerv\'{a}ry graph have a specific position with
respect to the maximum stable sets.

\begin{lemma}
\cite{levm2}\label{lem4} If $G$ is a K\"{o}nig-Egerv\'{a}ry graph, then for
any $S\in \Omega (G)$ every maximum matching of $G$ is contained in $%
(S,V(G)-S)$.
\end{lemma}

\begin{proposition}
\label{prop12}If $M$ is a maximum matching in a graph $G$ and $H$ is a
subgraph of $G$ such that $M=(M\cap E(H))\cup (M\cap E(G-H))$, then $\mu
(G)=\mu (H)+\mu (G-H)$.
\end{proposition}

\setlength {\parindent}{0.0cm}\textbf{Proof.} Clearly, $M\cap E(H)$ and $%
M\cap E(G-H)$ are matchings in $H$ and $G-H$, respectively. Let $M_{1},M_{2}$
be maximum matchings in $H$ and $G-H$, respectively. If $\mu (H)=\left|
M_{1}\right| >\left| M\cap E(H)\right| $, or $\mu (G-H)=\left| M_{2}\right|
>\left| M\cap E(G-H)\right| $, then $M_{1}\cup M_{2}$ is a matching in $G$
of cardinality larger than $\left| M\right| $, in contradiction with $\left|
M\right| =\mu (G)$. Therefore, $\mu (G)=\mu (H)+\mu (G-H)$. \rule{2mm}{2mm}%
\setlength
{\parindent}{3.45ex}

\begin{proposition}
\label{prop14}If $M$ is a maximum matching in a K\"{o}nig-Egerv\'{a}ry graph 
$G$, and $H$ is a subgraph of $G$ such that $M=(M\cap E(H))\cup (M\cap
E(G-H))$, then

($\mathit{i}$) $H$ and $G-H$ are K\"{o}nig-Egerv\'{a}ry graphs;

($\mathit{ii}$) $\alpha (G)=\alpha (H)+\alpha (G-H)$.
\end{proposition}

\setlength {\parindent}{0.0cm}\textbf{Proof.} Let $S\in \Omega
(G),S_{1}=S\cap V(H)$ and $S_{2}=S\cap V(G-H)$. By Lemma \ref{lem4}, $%
M\subseteq (S,V(G)-S)$, and according to Proposition \ref{prop11}($\mathit{ii%
}$), $G=H_{1}*H_{2}$, where $V(H_{1})=S\in \Omega (G)$ and $\left|
V(H_{1})\right| \geq \mu (G)=\left| M\right| =\left| V(H_{2})\right| $.
Hence, we infer that: $V(H)=S_{1}\cup (V(H_{2})-V(G-H)),S_{1}$ is stable in $%
H,M\cap E(H)$ is a matching in $H$ of size $\left| V(H_{2})-V(G-H)\right| $,
and $\left| S_{1}\right| \geq \left| V(H_{2})-V(G-H)\right| $, i.e., $H$ is
a K\"{o}nig-Egerv\'{a}ry graph, according to Proposition \ref{prop11}($%
\mathit{iii}$). Similarly, $G-H$ is also a K\"{o}nig-Egerv\'{a}ry graph.
Since, by Proposition \ref{prop12}, $\mu (G)=\mu (H)+\mu (G-H)$ and all $%
G,H,G-H$ are K\"{o}nig-Egerv\'{a}ry graphs, we may conclude that $\alpha
(G)=\alpha (H)+\alpha (G-H)$. \rule{2mm}{2mm}\setlength
{\parindent}{3.45ex}

\begin{lemma}
\label{lem1}If $H$ is a subgraph of $G$, such that $\alpha (G)=\alpha
(H)+\alpha (G-H)$ and $G$ is $\alpha ^{++}$-stable, then $H$ is $\alpha
^{++} $-stable, as well.
\end{lemma}

\setlength {\parindent}{0.0cm}\textbf{Proof.} Since $\alpha (G)=\alpha
(H)+\alpha (G-H),$ it follows that any $S\in \Omega (G)$ satisfies $\left|
S\cap V(H)\right| =\alpha (H)$. So, if $\alpha (H+e_{1}+e_{2})<\alpha (H)$,
for some $e_{1},e_{2}\in E(\overline{H})$, it follows that $\alpha
(G+e_{1}+e_{2})<\alpha (G)$, as well. \rule{2mm}{2mm}%
\setlength
{\parindent}{3.45ex}

\begin{lemma}
\label{lem2}If $G$ is of order $6$, has a Hamiltonian path and $\alpha (G)=3$%
, then $G$ is not $\alpha ^{++}$-stable.
\end{lemma}

\setlength {\parindent}{0.0cm}\textbf{Proof.} Suppose that $%
V(G)=\{v_{i}:1\leq i\leq 6\}$ and $(v_{i},v_{i+1})\in E(G)$ for any $i\in
\{1,...,5\}$. Then $H=G+v_{1}v_{3}+v_{4}v_{6}$ has $\alpha (H)=2$, i.e., $G$
is not $\alpha ^{++}$-stable. \rule{2mm}{2mm}\setlength
{\parindent}{3.45ex}

\begin{proposition}
\label{prop5}If $G\neq K_{n}-e,n=2,3$ is an $\alpha ^{++}$-stable
K\"{o}nig-Egerv\'{a}ry graph, then $G$ has a perfect matching consisting of
only pendant edges.
\end{proposition}

\setlength {\parindent}{0.0cm}\textbf{Proof.} If $G$ is $\alpha ^{++}$%
-stable then, clearly, $G$ is $\alpha _{P_{3}}^{+}$-stable too. Hence, by
Proposition \ref{prop10} if $G\neq K_{n}-e$ then it is also $\alpha ^{+}$%
-stable. It is not difficult to check that $K_{n}-e$ can be a
K\"{o}nig-Egerv\'{a}ry graph only for $n=2,3$. Therefore, if a
K\"{o}nig-Egerv\'{a}ry graph $G\neq K_{n}-e,n=2,3$ is $\alpha ^{++}$-stable,
then it is also $\alpha ^{+}$-stable. Now Proposition \ref{prop1} ensures
that $G$ has a perfect matching, say $M=\{a_{i}b_{i}:1\leq i\leq \alpha
(G)\} $. According to Proposition \ref{prop11} and Lemma \ref{lem4}, we may
assume that $S=\{a_{i}:1\leq i\leq \alpha (G)\}\in \Omega (G)$. We show that 
$M$ consists of only pendant edges. Suppose, on the contrary, that some $%
a_{k}b_{k}\in M$ is not pendant.\setlength
{\parindent}{3.45ex}

\begin{figure}[h]
\setlength{\unitlength}{1cm}%
\begin{picture}(5,2)\thicklines
  \multiput(2,0.5)(1,0){2}{\circle*{0.29}}
  \multiput(2,1.5)(1,0){2}{\circle*{0.29}}
  \put(2,0.5){\line(1,0){1}} 
  \put(2,0.5){\line(0,1){1}}
  \put(2,0.5){\line(1,1){1}}
  \put(3,0.5){\line(0,1){1}}
  \multiput(2,1.5)(0.125,0){8}{\circle*{0.07}}
  \multiput(2,1.5)(0.125,-0.125){8}{\circle*{0.07}}
  \put(1.5,1){\makebox(0,0){$(a)$}}
  \put(2,0){\makebox(0,0){$b_{i}$}}
  \put(2,1.9){\makebox(0,0){$a_{i}$}}
  \put(3,0){\makebox(0,0){$b_{k}$}}
  \put(3,1.9){\makebox(0,0){$a_{k}$}}

  \multiput(5,0.5)(1,0){2}{\circle*{0.29}}
  \multiput(5,1.5)(1,0){2}{\circle*{0.29}}
  \put(5,0.5){\line(0,1){1}} 
  \put(5,0.5){\line(1,1){1}}
  \put(5,1.5){\line(1,-1){1}}
  \put(6,0.5){\line(0,1){1}} 
  \multiput(5,0.5)(0.125,0){8}{\circle*{0.07}}
  \multiput(5,1.5)(0.125,0){8}{\circle*{0.07}}
  \put(4.5,1){\makebox(0,0){$(b)$}}
  \put(5,0){\makebox(0,0){$b_{i}$}}
  \put(5,1.9){\makebox(0,0){$a_{i}$}}
  \put(6,0){\makebox(0,0){$b_{k}$}}
  \put(6,1.9){\makebox(0,0){$a_{k}$}}

  \multiput(8,0.5)(1,0){3}{\circle*{0.29}}
  \multiput(8,1.5)(1,0){3}{\circle*{0.29}}
  \put(8,0.5){\line(0,1){1}} 
  \put(9,0.5){\line(0,1){1}}
  \put(10,0.5){\line(0,1){1}}
  \put(8,0.5){\line(1,1){1}} 
  \put(9,0.5){\line(1,1){1}}
  \multiput(8,1.5)(0.125,0){8}{\circle*{0.07}}
  \multiput(9,0.5)(0.125,0){8}{\circle*{0.07}}
  \put(7.5,1){\makebox(0,0){$(c)$}}
  \put(8,0){\makebox(0,0){$b_{i}$}}
  \put(8,1.9){\makebox(0,0){$a_{i}$}}
  \put(9,0){\makebox(0,0){$b_{k}$}}
  \put(9,1.9){\makebox(0,0){$a_{k}$}}
  \put(10,0){\makebox(0,0){$b_{j}$}}
  \put(10,1.9){\makebox(0,0){$a_{j}$}}
  
 \end{picture}
\caption{Non-$\alpha ^{++}$-stable K\"{o}nig-Egerv\'{a}ry graphs.}
\label{fig1323245}
\end{figure}
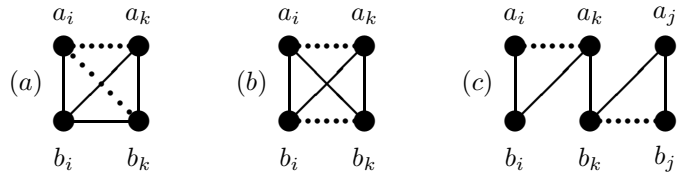

\textit{Case 1}. There exists some $b_{i}$ such that $a_{k}b_{i},b_{i}b_{k}%
\in E(G)$ (see Figure \ref{fig1323245}\textit{(a)}). If $H=G[%
\{a_{k},b_{i},a_{i},b_{k}\}]$, then Proposition \ref{prop14}($\mathit{ii}$)
implies that $\alpha (G)=\alpha (H)+\alpha (G-H)$. Since $H$ is not $\alpha
^{++}$-stable, it follows, by Lemma \ref{lem1}, that $G$ could not be $%
\alpha ^{++}$-stable, in contradiction with the premises on $G$.

\textit{Case 2}. There exist $a_{i}b_{i}\in M$ with $a_{k}b_{i},a_{i}b_{k}%
\in E(G)$ (see Figure \ref{fig1323245}\textit{(}$\mathit{b}$\textit{)}). If $%
H=G[\{a_{k},b_{i},a_{i},b_{k}\}]$, then Proposition \ref{prop14}($\mathit{ii}
$) ensures that $\alpha (G)=\alpha (H)+\alpha (G-H)$. Since $H$ is not $%
\alpha ^{++}$-stable, it follows, by Lemma \ref{lem1}, that $G$ could not be 
$\alpha ^{++}$-stable, in contradiction with the premises on $G$.

\textit{Case 3}. There exist $a_{i},b_{i},a_{j},b_{j}$, such that $%
a_{i}b_{i},a_{j}b_{j}\in M$ and $a_{k}b_{i},a_{j}b_{k}\in E(G)$. In
addition, we can assume that $b_{i}b_{k},b_{k}b_{j}\notin E(G)$, otherwise
we return to \textit{Case 1.} Hence, $H=G[%
\{a_{i},a_{k},a_{j},b_{i},b_{k},b_{j}\}]$ contains a path on $6$ vertices
(see Figure \ref{fig1323245}\textit{(}$\mathit{c}$\textit{)}). Since $\alpha
(H)=\left| \{a_{i},a_{k},a_{j}\}\right| =3$, Lemma \ref{lem2} implies that $%
H $ is not $\alpha ^{++}$-stable, and because $\alpha (G)=\alpha (H)+\alpha
(G-H)$ is true according to Proposition \ref{prop14}($\mathit{ii}$), we get,
by Lemma \ref{lem1}, that $G$ cannot be $\alpha ^{++}$-stable, in
contradiction with the premises on $G$.

Thus, $M$ must consist of only pendant edges. \rule{2mm}{2mm}%
\setlength
{\parindent}{3.45ex}\newline

\begin{figure}[h]
\setlength{\unitlength}{1cm}%
\begin{picture}(5,1.2)\thicklines
  \multiput(3,0)(1,0){4}{\circle*{0.29}}
  \multiput(3,1)(1,0){4}{\circle*{0.29}}
  \put(3,0){\line(1,0){3}}
  \put(3,1){\line(1,0){3}}
  \put(4,0){\line(0,1){1}} 
  \put(4,0){\line(1,1){1}} 
  \put(5,0){\line(0,1){1}}

\put(2.7,0){\makebox(0,0){$a$}}
\put(2.7,1){\makebox(0,0){$b$}}
\put(6.35,0){\makebox(0,0){$c$}}
\put(6.35,1){\makebox(0,0){$d$}}

  \multiput(8,0)(1,0){4}{\circle*{0.29}}
  \multiput(9,1)(1,0){2}{\circle*{0.29}}
  \put(8,0){\line(1,0){3}}
  \put(9,1){\line(1,0){1}}
  \put(9,1){\line(1,-1){1}} 
  \put(9,0){\line(0,1){1}}

\put(2.25,0.5){\makebox(0,0){$G_{1}$}}
\put(7.5,0.5){\makebox(0,0){$G_{2}$}}
 \end{picture}
\caption{K\"{o}nig-Egerv\'{a}ry graphs with a perfect matching consisting of
pendant edges.}
\label{fig80}
\end{figure}
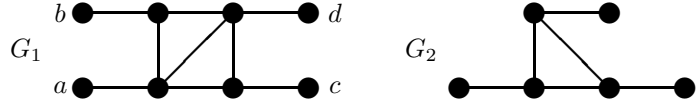

It is worth observing that : ($\mathit{a}$) Proposition \ref{prop5} fails
for non-K\"{o}nig-Egerv\'{a}ry graphs; e.g., $C_{5}$ is $\alpha ^{++}$%
-stable and has no perfect matching; ($\mathit{b}$) the converse of
Proposition \ref{prop5}, within the class of K\"{o}nig-Egerv\'{a}ry graphs,
is not generally true. For instance, the graph $G_{1}$ in Figure \ref{fig80}
has a perfect matching consisting of only pendant edges and it is not $%
\alpha ^{++}$-stable (because $\alpha (G_{1}+ad+bc)<\alpha (G_{1})$), while
the graph $G_{2}$ in the same figure has a perfect matching consisting of
only pendant edges, and it is also $\alpha ^{++}$-stable.

However, we can show that:

\begin{proposition}
\label{prop6}Any graph that has a perfect matching consisting of only
pendant edges is $\alpha _{P_{3}}^{+}$-stable.
\end{proposition}

\setlength {\parindent}{0.0cm}\textbf{Proof.} Let $G$ be a graph that has a
perfect matching $M=\{a_{i}b_{i}:1\leq i\leq \mu (G)\}$, consisting of only
pendant edges, and suppose that $S_{0}=\{a_{i}:1\leq i\leq \mu (G)\}\in
\Omega (G)$. Let denote $H=G+e_{1}+e_{2}$, where $e_{1},e_{2}\in E(\overline{%
G})$ are such that they have a common endpoint, say $e_{1}=uv,e_{2}=vw$. 
\setlength
{\parindent}{3.45ex}We distinguish between the following cases:

\textit{Case 1}. If $u,v,w\notin S_{0}$, then $S_{0}\in \Omega (H)$, and $%
\alpha (H)=\alpha (G)$.

\textit{Case 2}. If $u,v\notin S_{0}$ or $u,w\notin S_{0}$, then $S_{0}\in
\Omega (H)$, and $\alpha (H)=\alpha (G)$.

\textit{Case 3}. If $u\notin S_{0}$ and $v=a_{i},w=a_{j}$, then $S_{0}\cup
\{b_{j}\}-\{a_{j}\}\in \Omega (H)$, and $\alpha (H)=\alpha (G)$.

\textit{Case 4}. If $v\notin S_{0}$ and $u=a_{i},w=a_{j}$, then $S_{0}\in
\Omega (H)$, and $\alpha (H)=\alpha (G)$.

\textit{Case 5}. If $u=a_{i},v=a_{j},w=a_{k}$, then $S_{0}\cup
\{b_{j}\}-\{a_{j}\}\in \Omega (H)$, and $\alpha (H)=\alpha (G)$.

Consequently, $G$ is $\alpha _{P_{3}}^{+}$-stable. \rule{2mm}{2mm}

\begin{theorem}
\label{th1}A graph $G$ that has a perfect matching consisting of only
pendant edges is $\alpha ^{++}$-stable if and only if $G$ contains no cycle
on $4$ vertices.
\end{theorem}

\setlength {\parindent}{0.0cm}\textbf{Proof.} Let $M=\{a_{i}b_{i}:1\leq
i\leq \mu (G)\}$ be the perfect matching of $G$. Without loss of generality,
we can assume that $S_{0}=\{a_{i}:1\leq i\leq \mu (G)\}\in \Omega (G)$.%
\setlength
{\parindent}{3.45ex}

Suppose, on the contrary, that there is 
\[
D=\{b_{i},b_{j},b_{k},b_{m}\}\ with\
b_{i}b_{j},b_{j}b_{k},b_{k}b_{m},b_{m}b_{i}\in E(G), 
\]
i.e., $G[D]$ contains a Hamiltonian cycle. If $H=G[\{a_{i},a_{j},a_{k},a_{m}%
\}\cup D]$, then $\alpha (G)=\alpha (G-H)+\alpha (H)$, since any $S\in
\Omega (G)$ satisfies $\left| S\cap \{a_{q},b_{q}\}\right| =1$ for each $%
a_{q}b_{q}\in M$. On the one hand, by Lemma \ref{lem1}, $H$ should be $%
\alpha ^{++}$-stable. On the other hand, $\alpha
(H+a_{i}a_{k}+a_{j}a_{m})=3<\alpha (H)=4$, which brings a contradiction.
Therefore, $G$ has no cycle on $4$ vertices.

Conversely, let $G$ be such that no $4$ vertices span a cycle. Assume, on
the contrary, that $G$ is not $\alpha ^{++}$-stable, i.e., there are $%
e_{1},e_{2}\in E(\overline{G})$ such that $H=G+e_{1}+e_{2}$ has $\alpha
(H)<\alpha (G)$. If at least one of $e_{1},e_{2}$ joins two vertices from $%
\{b_{i}:1\leq i\leq \mu (G)\}$, or one from $\{a_{i}:1\leq i\leq \mu (G)\}$
and the second from $\{b_{i}:1\leq i\leq \mu (G)\}$, then $\alpha (H)=\alpha
(G)$. According to Proposition \ref{prop6}, the same result follows if $%
e_{1},e_{2}$ have a common endpoint. Suppose that $%
e_{1}=a_{i}a_{j},e_{2}=a_{k}a_{m}$ and $a_{i},a_{j},a_{k},a_{m}$ are
pairwise distinct. Hence, we get:

\begin{itemize}
\item  $b_{i}b_{k}\in E(G)$, otherwise any $S\in \Omega (G)$ containing $%
\{a_{j},a_{m},b_{i},b_{k}\}$ is stable in $H$;

\item  $b_{j}b_{m}\in E(G)$, otherwise any $S\in \Omega (G)$ containing $%
\{a_{i},a_{k},b_{j},b_{m}\}$ is stable in $H$;

\item  $b_{i}b_{m}\in E(G)$, otherwise any $S\in \Omega (G)$ containing $%
\{a_{j},a_{k},b_{i},b_{m}\}$ is stable in $H$;

\item  $b_{j}b_{k}\in E(G)$, otherwise any $S\in \Omega (G)$ containing $%
\{a_{i},a_{m},b_{j},b_{k}\}$ is stable in $H$.
\end{itemize}

It follows that $b_{i}b_{k},b_{j}b_{k},b_{j}b_{m},b_{i}b_{m}\in E(G)$, i.e., 
$\{b_{i},b_{j},b_{k},b_{m}\}$ spans a $4$-cycle in $G$, in contradiction
with the premises on $G$. Consequently, $G$ is $\alpha ^{++}$-stable. \rule%
{2mm}{2mm}\newline

Combining Proposition \ref{prop5} and Theorem \ref{th1}, we obtain the
following characterization of $\alpha ^{++}$-stable K\"{o}nig-Egerv\'{a}ry
graphs.

\begin{theorem}
\label{th3}A K\"{o}nig-Egerv\'{a}ry graph is $\alpha ^{++}$-stable if and
only if it has a perfect matching consisting of only pendant edges and
contains no cycle on $4$ vertices.
\end{theorem}

Recall that a graph $G$ is called: ($\mathit{a}$) \textit{well-covered} if
every maximal stable set of $G$ is also a maximum stable set, i.e., it is in 
$\Omega (G)$, \cite{plum}; ($\mathit{b}$) \textit{very well-covered}
provided $G$ is well-covered and $\left| V(G)\right| =2\alpha (G)$, \cite
{fav1}. The following result extends the characterization that Finbow,
Hartnell and Nowakowski give in \cite{finhart1} for well-covered graphs
having the girth $\geq 6$.

\begin{proposition}
\label{prop13}\label{prop9}Let $G$ be a graph of girth $\geq 6$, which is
isomorphic to neither $C_{7}$ nor $K_{1}$. Then the following assertions are
equivalent:

($\mathit{i}$) $G$ is well-covered;

($\mathit{ii}$) $G$ has a perfect matching consisting of pendant edges;

($\mathit{iii}$) $G$ is very well-covered;

($\mathit{iv}$) $G$ is a K\"{o}nig-Egerv\'{a}ry $\alpha _{0}^{+}$-stable
graph with exactly $\alpha (G)$ pendant vertices;

($\mathit{v}$) $G$ is a K\"{o}nig-Egerv\'{a}ry $\alpha ^{++}$-stable graph.
\end{proposition}

\setlength {\parindent}{0.0cm}\textbf{Proof.} The equivalences ($\mathit{i}$%
) $\Leftrightarrow $ ($\mathit{ii}$) $\Leftrightarrow $ ($\mathit{iii}$) are
done in \cite{finhart1}. In \cite{levm12} it has been proved that ($\mathit{%
iii}$) $\Leftrightarrow $ ($\mathit{iv}$). Finally, ($\mathit{ii}$) $%
\Leftrightarrow $ ($\mathit{v}$) is true by Theorem \ref{th3}. \rule%
{2mm}{2mm}\setlength {\parindent}{3.45ex}

\begin{corollary}
\label{cor1}For a bipartite graph $G$ the following assertions are
equivalent:

($\mathit{i}$) $G$ is $\alpha ^{++}$-stable;

($\mathit{ii}$) $G$ is $C_{4}$-free and has a perfect matching consisting of
only pendant edges;

($\mathit{iii}$) $G$ is $C_{4}$-free and well-covered.
\end{corollary}

\begin{corollary}
$C_{n}$ is $\alpha ^{++}$-stable if and only if $n$ is odd.
\end{corollary}

\setlength {\parindent}{0.0cm}\textbf{Proof.} For any $n\geq 2,C_{2n}$ is
not an $\alpha ^{++}$-stable graph according to Corollary \ref{cor1}. 
\setlength
{\parindent}{3.45ex}

Assume, on the contrary, that $C_{2n+1}$ is not an $\alpha ^{++}$-stable
graph. Hence, there are $e_{1},e_{2}\in E(\overline{C_{2n+1}}%
),e_{1}=xy,e_{2}=uv$ such that $\alpha (C_{2n+1}+e_{1}+e_{2})<\alpha
(C_{2n+1})$. We may suppose, without loss of generality, that $x=v_{1}\neq u$%
. Now, if $u=v_{2i+1}$ (for some $i\neq 0$), then $x,u\notin
S=\{v_{2i}:1\leq i\leq n\}\in \Omega (C_{2n+1})$, and if $u=v_{2i}$ (for
some $i\neq 0$), then $x,y\notin S=\{v_{2},v_{4},...,v_{2i-2}\}\cup
\{v_{2i+1},v_{2i+3},...,v_{2n+1}\}\in \Omega (C_{2n+1})$. Hence, we infer
that $S$ is stable in $C_{2n+1}+e_{1}+e_{2}$, as well, in contradiction with 
$\alpha (C_{2n+1}+e_{1}+e_{2})<\alpha (C_{2n+1})$. Therefore, $C_{2n+1}$ is $%
\alpha ^{++}$-stable. \rule{2mm}{2mm}\setlength
{\parindent}{3.45ex}\newline

Combining Corollary \ref{cor1} and Proposition \ref{prop9} we get the
following extension of one Ravindra's theorem, \cite{rav1}, where he proved
the first three equivalences.

\begin{corollary}
For a tree $T$ the following assertions are equivalent:

($\mathit{i}$) $T$ is well-covered;

($\mathit{ii}$) $T$ has a perfect matching consisting of pendant edges;

($\mathit{iii}$) $T$ is very well-covered;

($\mathit{iv}$) $T$ is $\alpha ^{++}$-stable.
\end{corollary}

\section{Conclusions}

In this paper we keep investigating graphs whose stability number is
invariant with respect to some natural operations on graphs. While in \cite
{levm1}, \cite{levm2} we were interested in measuring the influence of
adding one edge to a graph, here we define a class of graphs whose stability
number is unaffected by two edges addition.

Further we concentrate on K\"{o}nig-Egerv\'{a}ry graphs, which is one of the
most attractive generalizations of bipartite graphs. One the one hand,
Proposition \ref{prop13} claims that for girth $\geq 6$, $\alpha ^{++}$%
-stable K\"{o}nig-Egerv\'{a}ry graphs and well-covered graphs are the same.
On the other hand, Theorem \ref{th3} shows that an $\alpha ^{++}$-stable
K\"{o}nig-Egerv\'{a}ry graph contains no cycle on $4$ vertices. It leaves an
interesting open question concerning interconnections between well-covered
graphs and $\alpha ^{++}$-stable K\"{o}nig-Egerv\'{a}ry graphs of girth $3$
or $5$.


\begin{thebibliography}{99}
\bibitem{berge}  C. Berge, \emph{Graphs}, North-Holland, Amsterdam, 1985.

\bibitem{dem}  R. W. Deming, \emph{Independence numbers of graphs - an
extension of the K\"{o}nig-Egerv\'{a}ry theorem}, Discrete Mathematics 
\textbf{27} (1979) 23-33.

\bibitem{ding}  G. Ding, \emph{Stable sets versus independent sets},
Discrete Mathematics \textbf{117} (1993) 73-87.

\bibitem{eger}  E. Egerv\'{a}ry, \emph{On combinatorial properties of
matrices}, Matematikai Lapok \textbf{38} (1931) 16-28.

\bibitem{fav1}  O. Favaron, \emph{Very well-covered graphs}, Discrete
Mathematics \textbf{42} (1982) 177-187.

\bibitem{finhart1}  A. Finbow, B. Hartnell and R. J. Nowakowski, \emph{A
characterization of well-covered graphs of girth 5 or greater}, Journal of
Combinatorial Theory Ser. B \textbf{57} (1993) 44-68.

\bibitem{gun}  G. Gunther, B. Hartnell, and D. F. Rall, \emph{Graphs whose
vertex independence number is unaffected by single edge addition or deletion}%
, Discrete Applied Mathematics \textbf{46} (1993) 167-172.

\bibitem{hayn}  T. W. Haynes, L. M. Lawson, R. C. Brigham and R. D. Dutton, 
\emph{Changing and unchanging of the graphical invariants: minimum and
maximum degree, maximum clique size, node independence number and edge
independence number}, Congressus Numerantium \textbf{72} (1990) 239-252.

\bibitem{koen}  D. K\"{o}nig, \emph{Graphen und Matrizen}, Matematikai Lapok 
\textbf{38} (1931) 116-119.

\bibitem{levm1}  V. E. Levit and E. Mandrescu, \emph{On }$\alpha $\emph{%
-stable graphs}, Congressus Numerantium \textbf{124} (1997) 33-46.

\bibitem{levm12}  V. E. Levit and E. Mandrescu, \emph{Well-covered and
K\"{o}nig-Egerv\'{a}ry graphs}, Congressus Numerantium \textbf{130} (1998)
209-218.

\bibitem{levm2}  V. E. Levit and E. Mandrescu, \emph{On }$\alpha ^{+}$\emph{%
-stable K\"{o}nig-Egerv\'{a}ry graphs}, The Ninth SIAM\ Conference on
Discrete Mathematics, University of Toronto, Canada (1998), Los Alamos
Archive, prE-print math.CO/9912022, 1999, 13 pp.

\bibitem{plum}  M. D. Plummer, \emph{Some covering concepts in graphs},
Journal of Combinatorial Theory \textbf{8} (1970) 91-98.

\bibitem{rav1}  G. Ravindra, \emph{Well-covered graphs}, Journal of
Combinatorial Information System Sciences \textbf{2} (1977) 20-21.

\bibitem{ster}  F. Sterboul, \emph{A characterization of the graphs in which
the transversal number equals the matching number}, Journal of Combinatorial
Theory B \textbf{27} (1979) 228-229.
\end{thebibliography}
\end{document}